\documentclass[12pt,a4paper]{article}

\makeatletter
\usepackage[usenames,dvips]{color}

\let\L@TeX@startsection\@startsection
\renewcommand{\@startsection}[6]{\L@TeX@startsection{#1}{#2}{#3}{#4}{#5}%
{\color{red}#6}}
\let\L@TeXmaketitle\maketitle
\renewcommand{\maketitle}{{\color{Brown}\L@TeXmaketitle}}
\let\L@TeXverbatim@font\verbatim@font
\def\verbatim@font{\color{Sepia}\L@TeXverbatim@font}
\newcommand{\ajr@footercolour}{\color{green}}
\newcommand{\ajr@headercolour}{\color{red}}
\newcommand{\ajr@theoremcolour}{\color{blue}}

\renewcommand{\vec}[1]{\mbox{\boldmath$#1$}}
\newcommand{\Ord}[1]{{\cal O}\left(#1\right)}

\newtheorem{theorem}{\ajr@theoremcolour Theorem}

\newenvironment{proof}%
{\par\textbf{\ajr@theoremcolour Proof:}}%
{\hfill{\ajr@theoremcolour$\spadesuit$}\par}

\let\ajr@title\title
\renewcommand{\title}[1]{\ajr@title{\sf #1}}

\let\ajr@marginpar\marginpar
\renewcommand{\marginpar}[1]{\ajr@marginpar{\raggedright\footnotesize\sf #1}}

\renewcommand{\MakeUppercase}[1]{\ajr@headercolour\textsf{#1}}
\def\@oddfoot{\hfill\footnotesize\sf\ajr@footercolour AJ Roberts, \today}

\let\ajr@@openbib@code\@openbib@code
\renewcommand{\@openbib@code}{%
    \ajr@@openbib@code
    \addcontentsline{toc}{section}{\refname}
    }

\usepackage{url}
\makeatother

\newcommand{\RR}{{\mathrm{I\!\!R}}}

\begin{document}
\title{A flexible error estimate for the application of centre manifold
theory}
\author{Zhenquan Li \and A.J.  Roberts\thanks{Dept of Mathematics \&
Computing, University of Southern Queensland, Toowoomba, Queensland
4350, \textsc{Australia};  \protect\url{mailto:zhen.li@auckland.ac.nz} and
\protect\url{mailto:aroberts@usq.edu.au} respectively.} }
\maketitle

\begin{abstract}
In applications of centre manifold theory we need more flexible error
estimates than that provided by, for example, the Approximation
Theorem~3 by Carr~\cite{Carr81, Carr83b}.  Here we extend the theory
to cover the case where the order of approximation in parameters and
that in dynamical variables may be completely different.  This allows,
for example, the effective evaluation of low-dimensional dynamical
models at finite parameter values.
\end{abstract}

\paragraph{Maths Subj.\ Class:} 37L10, 37N10.

\tableofcontents

\section{Introduction}

Interest in the dynamical behaviour of a physical system usually lies
in the relatively low-dimensional evolution after heavily damped modes
have become insignificant.
There are many successful applications of centre manifold techniques
to create models of these relatively simple dynamics.
We here mention some applications to physical fluid mechanics.
Iooss~\cite{Iooss92, Iooss90}, and Laure~\cite{Laure88} analysed the
dynamics of Taylor vortices in Taylor-Couette flow, whereas
Chossat~\cite{Chossat87} and Hill~\cite{Hill91} discuss the
non-axisymmetric dynamics involving mode competition by using centre
manifold theory.
Mode interactions in the dynamics of convection in porous media are
analysed with centre manifolds by Neel~\cite{Neel90, Neel91, Neel93}
and Graham \& Steen~\cite{Graham92}.
Arneodo {\em et al}~\cite{Arneodo85c, Arneodo85b, Arneodo85a} reduced
the dynamics of triple convection down to a set of three coupled
\textsc{ode}s, numerically verified the modelling and then proved the
existence of chaos.
Roberts \emph{et al} discussed centre manifolds of forced dynamical
systems \cite{Cox91}, and derived low-dimensional models using centre
manifold techniques for contaminant dispersion in channels
\cite{Mercer90}, shear dispersion in pipes \cite{Mercer94a}, thin film
fluid dynamics \cite{Roberts96b}, coating flows over a curved
substrate in space \cite{Roy96, Zhenquan98}, and Mei, Roberts \&
Li~\cite{Mei94, Roberts99c} derived models for turbulent shallow
water flow written in terms of vertically averaged quantities derived
from the $k$-$\epsilon$ model for turbulent flow.
Such applications assure us that centre manifold theory provides a
useful route to the low-dimensional modelling of high-dimensional
dynamical systems.

Centre manifold theory guarantees the existence of low-dimensional
models, matches the solutions of original and the low-dimensional
systems, and quantifies errors in the approximation.
Algebraic techniques to construct low-dimensional models are based
upon the theory.
In problems specified in the standard form~(\ref{genpr}), the centre
manifold may be calculated simply by iteration, see Carr~\cite{Carr81}
for example.
For the more directly applicable form~$\dot{\vec{u}}={\cal
L}\vec{u}+\vec{f}(\vec{u},\vec\epsilon)$, solutions in the form of an
asymptotic power series are found using methods developed by Coullet
\& Spiegel~\cite{Coullet83} (and reinvented by Leen~\cite{Leen93}).
The derivation of initial conditions for such low-dimensional models
is given through projecting the initial condition of the system onto
the centre manifold \cite{Cox93b, Roberts89, Roberts92, Roberts97b}.
But many of the applications require more flexible errors estimates.
For example, physical models recoverd by evaluation at a finite value
of a supposedly asymptotically small parameter often need high order
approximations in the parameter \cite{Mercer90, Roberts94c, Mei94,
Roberts96b, Roberts99c, Roberts99b}.
Thus asymptotic errors estimates need to be made to high order in some
parameters and only low order in other variables.
In Section~\ref{exth} we extend a theorem of Carr \&
Muncaster~\cite{Carr81,Carr83a} to rigorously support such flexible
approximations.

\section{A simple example}

To introduce the issues, consider the well known prototype bifurcation
problem
\begin{equation}
\dot{x}=\epsilon x-xy\,,\quad
\dot{y}=-y+x^2\,,\label{eq:prot}
\end{equation}
where $\epsilon$ is a parameter.  By adjoining the trivial equation,
it becomes:
\begin{equation}
\dot{\epsilon}=0\,,\quad
\dot{x}=\epsilon x-xy\,,\quad
\dot{y}=-y+x^2\,.\label{eq:prote}
\end{equation}
According to the distribution of the eigenvalues of the system, we may
seek a centre manifold of form $y=h(x,\epsilon)$.
Substituting this into the system~(\ref{eq:prote}) we deduce that $h$
must satisfy
\begin{equation}
h=x^2-\frac{\partial h}{\partial x}x(\epsilon-h)\,.
\label{eq:egres}
\end{equation}
Solving this iteratively leads to the approximations
\begin{equation}
h^{(0)}=0\,,\quad
h^{(1)}=x^2\,,\quad
h^{(2)}=x^2-2\epsilon x^2+2x^4\,,\quad\mbox{etc.}
\label{eq:egcm}
\end{equation}
Now elementary calculation shows that the above approximations
$h^{(n)}$ satisfy~(\ref{eq:egres}) to a residual~$\Ord{
|(\epsilon,x)|^{n+2}}$ as $(\epsilon,x) \rightarrow {\vec 0}$; an
error equivalently expressed as $\Ord{\epsilon^{n+2}+x^{n+2}}$ since a
term $c\epsilon^{p'} x^{q'}$ (for some constant $c\neq 0$) is
$\Ord{\epsilon^p+x^q}$ only if $p'/p+q'/q\geq 1$.
Therefore the centre manifold is $y=h(x,\epsilon)
=h^{(n)}+\Ord{|(\epsilon,x)|^{n+2}}$ by, for example,
\textsc{Theorem}~3 of \cite[p264]{Carr83a}.

The limitation in applications is that the established theorem on
approximation strongly couples the order of truncation in both
parameters and variables---the ``weight'' of the parameter and
variable is the same in the error estimate.
Some flexibilty may be introduced by a nonlinear transformation of the
parameters; for example, introduce $\delta=\sqrt\epsilon$ and instead
of~(\ref{eq:prote}) study
\begin{equation}
\dot{\delta}=0\,,\quad
\dot{x}=\delta^2 x-xy\,,\quad
\dot{y}=-y+x^2\,.\label{eq:protd}
\end{equation}
The resultant iterative solution of~(\ref{eq:egres}) is identical
to~(\ref{eq:egcm}).
The only difference is that the approximation theorem asserts the
errors in $h^{(n)}$ are $\Ord{|(\delta,x)|^{2n+2}}$ as $(\delta,x)
\rightarrow {\vec 0}$, that is, $\Ord{\epsilon^{n+1}+x^{2n+2}}$.
Thus certain trivial nonlinear transformations make no difference to
the algebraic analysis but do affect the error estimate.
There must be more flexibility in the errors than has so far been
proved.

A more flexible error bound to include the effects of both ${
\epsilon}$ and ${ x}$ to any desired orders is to express and seek
errors as $\Ord{x^q,\epsilon^p}$.
Note that $f=\Ord{x^q,\epsilon^p}$ means that any terms in $f$ of
the form $c\epsilon^{p'} x^{q'}$ (for some constant $c\neq 0$) must
satisfy $p'\geq p$ or $q'\geq q$.
For example, we may deduce, supported by Theorem~\ref{theorem:4}
herein,
\begin{equation}
    h=x^2(1-2\epsilon+4\epsilon^2)+x^4(2-16\epsilon+88\epsilon^2)
    +\Ord{x^6,\epsilon^3}\,.
    \label{eq:flex}
\end{equation}
This kind of error allows us separately to choose the orders of the
parameter~$\epsilon$ and the variable~$x$ which we want to include in
the centre manifold.
For example, here we may compute to higher orders in $\epsilon$,
observe the pattern of coefficients in this simple problem, and
realise \cite{Roberts85b} that the above is the low-order Taylor
polynomial in $\epsilon$ of
\begin{equation}
    h=\frac{x^2}{1+2\epsilon}
    +\frac{2x^4}{(1+2\epsilon)^2(1+4\epsilon)}
    +\Ord{x^6}\,.
    \label{eq:esum}
\end{equation}
Such approimations to high-order in parameters and low-order in
dynamical variables were used, before proof, and are essential to the
analyses in \cite{Mercer90, Roberts94c, Mei94, Roberts96b, Roberts99c,
Roberts99b}.

\section{The flexible extension}
\label{exth}

Consider dynamical systems expressed in the form
\begin{eqnarray}
\dot{\vec{x}}&=&A\vec{x}+\vec{f}(\vec{x},\vec{y},{\vec
\epsilon})\,,\nonumber\\
\dot{\vec{y}}&=&B\vec{y}+\vec{g}(\vec{x},\vec{y},{\vec
\epsilon})\,,\label{genpr}
\end{eqnarray}
where ${\vec x}\in\RR^m$, ${\vec y}\in\RR^n$, ${\vec
\epsilon}\in\RR^l$ and $A$ and $B$ are constant matrices such that all
the eigenvalues of $A$ have zero real parts, and all eigenvalues of
$B$ have negative real parts.  Functions ${\vec f}$ and ${\vec g}$ are
nonlinear for ${\vec x}$, ${\vec y}$, ${\vec \epsilon}$ and ${\vec
f}({\vec 0},{\vec 0},{\vec 0})={\vec g}({\vec 0},{\vec 0},{\vec 0})
={\vec f}'({\vec 0},{\vec 0},{\vec 0})={\vec g}'({\vec 0},{\vec
0},{\vec 0})={\vec 0}$ (where ${\vec f}'=[{\vec f}_{\vec x},{\vec
f}_{\vec y},{\vec f}_{\vec \epsilon}]$, and similarly for ${\vec g}'$
and other Jacobians).

For any function ${\vec \phi}:\RR^m\times\RR^l\rightarrow\RR^n$ which
is a continuously differentiable function and $\vec\phi(\vec 0,\vec 0)
={\vec \phi}'({\vec 0},{\vec 0}) ={\vec 0}$, define
\begin{displaymath}
({H}{\vec \phi})={\vec \phi}_{\vec x}({\vec x},{\vec
\epsilon})[A{\vec x}+{\vec f}({\vec x},{\vec \phi}({\vec x},{\vec
\epsilon}),{\vec \epsilon})]-B{\vec \phi}({\vec x},{\vec
\epsilon})-{\vec g}({\vec x},{\vec \phi}({\vec x},{\vec
\epsilon}),{\vec \epsilon})\,.
\end{displaymath}
Also let $\Ord{s^q,\epsilon^p}$ denote any terms of the form
$c\epsilon_1^{p_1}\ldots\epsilon_l^{p_l} s_1^{q_1}\ldots s_m^{q_m}$
(where the constant $c\neq 0$) which satisfy $p_1+\cdots+p_l\geq
p$ or $q_1+\cdots+q_m\geq q$ and $p_i,q_j\geq 0$.

\begin{theorem}[Approximation]\label{theorem:4}
Suppose that
\begin{displaymath}
({H}{\vec \phi})({\vec x},{\vec \epsilon})={\cal
O}(x^q,\epsilon^p) \quad\mbox{as}\quad ({\vec x},{\vec
\epsilon})\rightarrow {\vec 0} \,,
\end{displaymath}
where $p\geq 1$, $q>1$, then
\begin{displaymath}
|{\vec h}({\vec x},{\vec \epsilon})-{\vec \phi}({\vec x},{\vec
\epsilon})|=\Ord{x^q,\epsilon^p} \quad\mbox{as}\quad ({\vec
x},{\vec \epsilon})\rightarrow {\vec 0} \,.
\end{displaymath}
That is, the errors in the approximation $\vec\phi$ to a centre
manifold $\vec h$ is the same as the order of the residuals of the
equations of the dynamical system.
\end{theorem}

\begin{proof} This proof is adapted from Carr \cite[pp25--28]{Carr81}.

Let ${\vec \theta}:\RR^m\times\RR^l\rightarrow\RR^n$ be a continuously
differentiable function with compact support such that
\begin{displaymath}
{\vec \theta}({\vec x},{\vec \epsilon})={\vec \phi}({\vec x},{\vec
\epsilon})\quad\mbox{for}\quad|({\vec x},{\vec \epsilon})|
\quad\mbox{small}\,.
\end{displaymath}
Set
\begin{eqnarray}
{\vec N}({\vec x},{\vec \epsilon})&=&{\vec \theta}_{\vec x}({\vec
x},{\vec \epsilon}) \left[A{\vec x}+{\vec F}\left({\vec x},{\vec
\theta}({\vec x},{\vec \epsilon}),{\vec \epsilon}\right)\right]
\nonumber\\&&{}
-B{\vec \theta}({\vec x},{\vec \epsilon})-{\vec G}\left({\vec x},{\vec
\theta}({\vec x},{\vec \epsilon}),{\vec \epsilon}\right)\,,
\label{appphi}
\end{eqnarray}
where
\begin{eqnarray*}
{\vec F}({\vec x},{\vec y},{\vec \epsilon})&=&{\vec f}\left({\vec
x}\psi\left(\frac{\vec x}{\delta}\right),{\vec y},{\vec
\epsilon}\right)\,,
\\
{\vec G}({\vec x},{\vec y},{\vec
\epsilon})&=&{\vec g}\left({\vec x}\psi\left(\frac{\vec
x}{\delta}\right),{\vec y},{\vec \epsilon}\right)\,,
\end{eqnarray*}
where $\psi:\RR^m\rightarrow[0,1]$ is a infinitely differentiable
function with $\psi({\vec x})$=1 when $|{\vec x}|\leq 1$ and
$\psi({\vec x})=0$ when $|{\vec x}|\geq 2$ and $\delta$ is a positive
real number.
The properties of ${\vec F}$ and ${\vec G}$ are the same as in
\cite[p18]{Carr81} for $\vec\epsilon$ small.
So ${\vec N}({\vec x},{\vec \epsilon})=\Ord{x^q,\epsilon^p}$ as
$({\vec x},{\vec \epsilon})\rightarrow 0$.

For $a>0$ and $b>0$ let $\Gamma$ be the set of Lipschitz functions
${\vec h}:\RR^m\times\RR^l\rightarrow\RR^n$ with Lipschitz
constant $b$, $|{\vec h}({\vec x},{\vec \epsilon})|\leq a$ for $({\vec
x},{\vec \epsilon})\in \RR^m\times \RR^l$ and ${\vec h}({\vec 0},{\vec
0})={\vec 0}$.  With the supremum norm $\parallel.\parallel$, $\Gamma$
is a complete space.

For ${\vec h}\in\Gamma$ and ${\vec x}_0\in\RR^m$, let ${\vec
x}(t,{\vec x}_0,{\vec h})$ be the solution of
\begin{displaymath}
\dot{\vec x}=A{\vec x}+{\vec F}({\vec x},{\vec h},{\vec \epsilon})\,,
\quad{\vec x}(0,{\vec x}_0,{\vec h})={\vec x}_0\,.
\end{displaymath}
The bounds on ${\vec F}$ and ${\vec h}$ ensure that the solutions of
the above equation exists for all time $t$.  Define an operator $T$ on
$\Gamma$ by
\begin{displaymath}
(T{\vec h})({\vec x}_0)=\int_{-\infty}^0 e^{-Bs}{\vec G}\left({\vec
x}(s,{\vec x}_0,{\vec h}),{\vec h}({\vec x}(s,{\vec x}_0,{\vec
h})),{\vec \epsilon}\right)ds\,.
\end{displaymath}
We know that $T$ is a contraction mapping and the centre manifold
${\vec y}={\vec h}({\vec x},{\vec \epsilon})$ is a fixed point of $T$
for $a$, $b$ and $\delta$ small enough from the proof of the existence
theorem \cite[pp.16--19]{Carr81}.  Define
\begin{eqnarray*}
S{\vec Z}=T({\vec Z}+{\vec \theta})\,,
\end{eqnarray*}
for ${\vec Z}$ such that ${\vec Z}+{\vec \theta}\in \Gamma$.  The
domain of $S$ is a closed subset of $\Gamma$ since ${\vec
\theta}\in\Gamma$.  Since $|S{\vec Z}_1-S{\vec Z}_2|=|T({\vec
Z}_1+{\vec \theta})-T({\vec Z}_2+{\vec \theta})|$, thus $S$ is also a
contraction mapping.  For $K>0$ let\footnote{$K(x^q,\epsilon^p)$
denotes the product of $K$ and sum of finite terms
$c\epsilon_1^{p_1}\cdots\epsilon_l^{p_l} x_1^{q_1}\cdots x_m^{q_m}$,
where $p_1+\cdots+p_l\geq p$ or $q_1+\cdots+q_m\geq q$ and
$p_i$,$q_j\geq 0$, $c$ is a non-zero constant.}
\begin{eqnarray*}
Y=\left\{{\vec Z}\in\Gamma\mid|{\vec Z}({\vec x},{\vec \epsilon})|\le
K(x^q,\epsilon^p),\forall({\vec x},{\vec \epsilon})\in{\vec O}\subset
\RR^m\times\RR^l\right\}\,,
\end{eqnarray*}
where ${\vec O}$ is a neighbourhood of the origin in
$\RR^m\times\RR^l$.  Since ${\vec N}({\vec x},{\vec \epsilon})={\cal
O}(x^q,\epsilon^p)$ as $({\vec x},{\vec \epsilon})\rightarrow 0$, then
\begin{eqnarray}
|{\vec N}({\vec x},{\vec \epsilon})|\le
C_1(x^q,\epsilon^p)\,,\quad({\vec x},{\vec \epsilon})\in{\vec O}
\label{Nc1}
\end{eqnarray}
where $C_1$ is a constant.
Thus $Y$ is not empty because ${\vec N}\in Y\subset \Gamma$ by
defining $\vec{\theta}(\vec{x}, \vec{\epsilon})$ such that $C_1\le K$.
If we can find a constant $K$ such that $S$ maps $Y$ into $Y$, then
$\exists\ {\vec Z}_0\in Y$ is a fixed point of $S$, and
\begin{eqnarray*}
{\vec Z}_0=S({\vec Z}_0)=T({\vec Z}_0+{\vec \theta})-{\vec \theta}\,,
\\\mbox{that is,}\quad
T({\vec Z}_0+{\vec \theta})={\vec Z}_0+{\vec \theta}\,,
\end{eqnarray*}
i.e., ${\vec Z}_0+{\vec \theta}$ is a centre manifold
of~(\ref{genpr}), let ${\vec h}={\vec Z}_0+{\vec \theta}$,
\begin{eqnarray*}
|{\vec h}({\vec x},{\vec \epsilon})-{\vec \theta}({\vec x},{\vec
\epsilon})|={\vec Z}_0\le K(x^q,\epsilon^p)\,.
\end{eqnarray*}
To finish the proof define
\begin{eqnarray*}
{\vec Q}({\vec x},{\vec Z},{\vec \epsilon})&=&{\vec \theta}_{\vec
x}({\vec x},{\vec \epsilon}) \left[{\vec F}\left({\vec x},{\vec
\theta}+{\vec Z},{\vec \epsilon}\right)-{\vec F}\left({\vec x},{\vec
\theta},{\vec \epsilon}\right)\right] -{\vec N}({\vec x},{\vec
\epsilon})\\
& &{}+{\vec G}\left({\vec x},{\vec \theta}+{\vec Z},{\vec
\epsilon}\right)-{\vec G}\left({\vec x},{\vec \theta},{\vec
\epsilon}\right)\,.
\end{eqnarray*}
Then
\begin{eqnarray}
|{\vec Q}({\vec x},{\vec Z},{\vec \epsilon})|&\le&|{\vec Q}({\vec
x},{\vec 0},{\vec \epsilon})|+|{\vec Q}({\vec x},{\vec Z},{\vec
\epsilon})-{\vec Q}({\vec x},{\vec 0},{\vec \epsilon})|
\nonumber\\&&
=|{\vec N}({\vec x},{\vec \epsilon})|+|{\vec Q}({\vec x},{\vec
Z},{\vec \epsilon})-{\vec Q}({\vec x},{\vec 0},{\vec \epsilon})|\,.
\label{Qeq}
\end{eqnarray}
From the properties of $\vec F$ and $\vec G$ on p18 in \cite{Carr81}
and ${\vec \theta}'({\vec 0},{\vec 0})={\vec 0}$, we have
\begin{eqnarray}
|{\vec Q}({\vec x},{\vec Z},{\vec \epsilon})-{\vec Q}({\vec x},{\vec
0},{\vec \epsilon})|\le k(\delta)|{\vec Z}|\quad\mbox{for}\quad|{\vec
Z}|\le\delta\,.
\label{ineq2}
\end{eqnarray}
Using~(\ref{Nc1}), (\ref{Qeq}) and~(\ref{ineq2}),
\begin{eqnarray}
|{\vec Q}({\vec x},{\vec Z},{\vec \epsilon})|\le
(C_1+Kk(\delta))(x^q,\epsilon^p)\,,\quad\mbox{for}\quad{\vec Z}\in
Y\,.
\label{ineq3}
\end{eqnarray}
Using the same calculations as (2.5.9) on p27 \cite{Carr81}, for each
$r>0$, there is a constant $M(r)$ such that
\begin{eqnarray}
|{\vec x}(t,{\vec x}_0,{\vec \epsilon})|\le M(r)|{\vec x}_0|e^{-\gamma
t}\,,\quad t\le 0
\label{ineq4}
\end{eqnarray}
where $\gamma=r+2M(r)k(\delta)$ and ${\vec x}(t,{\vec x}_0,{\vec
\epsilon})$ is the solution of
\begin{eqnarray*}
\dot{\vec x}=A{\vec x}+{\vec F}({\vec x},{\vec Z}({\vec x},{\vec
\epsilon})+{\vec \theta}({\vec x},{\vec \epsilon}),{\vec
\epsilon})\,,\quad{\vec x}(0,{\vec x}_0,{\vec \epsilon})={\vec x}_0\,.
\end{eqnarray*}
Using (2.3.6) on p18 and (2.5.3) on p26 in \cite{Carr81},
and~(\ref{ineq3}), (\ref{ineq4}), if ${\vec Z}\in Y$
\begin{eqnarray*}
|(S{\vec Z})({\vec x}_0,{\vec
\epsilon})&\le&\left|\int_{-\infty}^0e^{-Bs}
(C_1+Kk(\delta))(x^q,\epsilon^p)\,ds\right|\\
&\le&\left|\int_{-\infty}^0e^{-Bs} (C_1+Kk(\delta))M(r)^q e^{-q\gamma
s}(x_0^q,\epsilon^p)\,ds\right|\\
&\le&C(C_1+Kk(\delta))M(r)^q(\beta-\gamma q)^{-1}(x_0^q,\epsilon^p)\,.
\end{eqnarray*}
provided $\delta$ and $r$ small enough so that $\beta-\gamma q>0$.
Choose ${\vec O}$ and $\delta$ small enough and $K$ large enough such
that $K\geq C(C_1+Kk(\delta))M(r)^q(\beta-\gamma q)^{-1}$.  Therefore
$S:Y\rightarrow Y$.  Hence Theorem~\ref{theorem:4} holds.
\end{proof}

More general theorems allowing varying orders of truncations within
parameters and dynamical variables may be also useful.
However, most such cases can be easily established by simple nonlinear
transformations of the parameters along the same lines as the example
in~(\ref{eq:protd}).

In applications, such as many fluid dynamics problems, we need theory
not only dealing with infinite dimensional problems, but also infinite
dimensional centre manifolds.
Carr \cite{Carr83b} presented the corresponding results for infinite
dimensional problems, and analysed two problems arising from partial
differential equations for finite dimensional centre manifolds.
The restrictions upon the nonlinear terms $\vec{f}$ is that $\vec{f}$
has order 2 continuous derivative and $\vec{f}(\vec{0})
=\vec{f}'(\vec{0}) =\vec{0}$.
More recently, Gallay~\cite{Gallay93} gave an extension of the
existence theorem to infinite dimensional centre manifolds, but a
bounded restriction on the nonlinear terms is required.
This condition limits its rigorous application.
Scarpellini \cite{Scarpellini91} apparently places significantly less
restrictions upon the nonlinearities in the dynamical equations, but
while he addresses infinite dimensional centre manifolds, the results
are severely constrained by requiring finite dimensional stable
dynamics.
H\u{a}r\u{a}gu\c{s} \cite{Haragus95, Haragus95b} has developed theory
supporting infinite dimensional models, such as the Korteweg-de Vries
equation, but only by placing extreme restrictions upon the linear
operators.
We identify the extension of the theorems to infinite dimensional
centre manifolds as a significant problem for future research.

\paragraph{Acknowledgement:} we thank the University of Southern
Queensland and the Australian Research Council for supporting this
research.

\bibliographystyle{plain}\bibliography{bib,ajr,new}

\end{document}